\nonstopmode

\documentclass[11pt]{article}

\usepackage[fleqn]{amsmath}

\usepackage{hyperref}
\hypersetup{
    colorlinks=true,
    linkcolor=blue,
    filecolor=magenta,      
    urlcolor=blue,
}

\usepackage{tcolorbox}   
\definecolor{mycolor}{rgb}{0.122, 0.435, 0.698}

  \usepackage{geometry}
\usepackage{amsmath,amssymb,amsthm}
\usepackage{graphics,epsfig,calc}

\usepackage{latexsym,epsfig,bm,amssymb}
\usepackage{xcolor}
\usepackage{amsthm,mathrsfs}

\usepackage{mathptmx}

\DeclareSymbolFont{AMSb}{U}{msb}m{n}
\DeclareSymbolFontAlphabet{\mathbb}{AMSb}

\newcommand{\beqn}{\begin{eqnarray}}
\newcommand{\eeqn}{\end{eqnarray}}
\newcommand{\be}{\begin{equation}}
\newcommand{\ee}{\end{equation}}

\newcommand{\bsp}{\begin{split}}
\newcommand{\esp}{\end{split}}

\newcommand{\ba}{\begin{array}}
\newcommand{\ea}{\end{array}}

\newcommand{\bpr}{\begin{proof}}
\newcommand{\epr}{\end{proof}}

\newcommand{\bH}{{\bf H}}

\newcommand{\bA}{{\bf A}}

\newcommand{\cD}{{\cal D}}
\newcommand{\bD}{{\bf D}}

\newcommand{\bF}{{\bf F}}

\newcommand{\cH}{{\cal H}}

\newcommand{\cK}{{\cal K}}

\newcommand{\cM}{{\cal M}}\newcommand{\bM}{{\bf M}}

\newcommand{\cO}{{\cal O}}

\newcommand{\bQ}{{\bf Q}}

\newcommand{\bS}{{\bf S}}
\newcommand{\rS}{{\rm S}}

\newcommand{\aX}{{\mathbb X}}

\newcommand{\aM}{{\mathbb M}}
\newcommand{\aQ}{{\mathbb Q}}

\newcommand{\bX}{{\bf X}}

\newcommand{\bY}{{\bf Y}}
\newcommand{\aY}{{\mathbb Y}}
\newcommand{\bZ}{{\bf Z}}

\newcommand{\al}{\alpha}

\newcommand{\ci}{\cite}
\newcommand{\de}{\delta}
\newcommand{{\De}}{{\Delta}}

\newcommand{\ds}{\displaystyle}
\newcommand{\fr}{\frac}

\newcommand{\ga}{\gamma}

\newcommand{\la}{\label}
\newcommand{\lam}{\lambda}

\newcommand{  \om}{  \omega}
\newcommand{  \Om}{  \Omega}

\newcommand{\vp}{\varphi}

\newcommand{ \ov}{ \overline}

\newcommand{\pa}{\partial}
\newcommand{\re}{\ref}

\newcommand{\we}{\wedge}

\newcommand{\Si}{\Sigma}

\newcommand{\ti}{\tilde}

\newcommand{\ve}{\varepsilon}

\newcommand\C{{\mathbb C}}
\newcommand\R{{\mathbb R}}

\newcommand{\bB}{{\bf B}}

\newcommand{\vka}{\varkappa}

\newcommand{\cm}{{\rm m}}

\newcommand{\bv}{{\bf v}}

\newcommand{\5}{{\hspace{0.5mm}}}

\newcommand{\const}{\mathop{\rm const}\nolimits}

\newcommand{\rRe}{{\rm Re\5}}
\newcommand{\rIm}{{\rm Im\5}}

\newtheorem{theorem}{Theorem}[section]

\renewcommand{\thetheorem}{\arabic{section}.\arabic{theorem}}

\newtheorem{definition}[theorem]{Definition}

\newtheorem{lemma}[theorem]{Lemma}
\newtheorem{example}[theorem]{Example}

\newtheorem{remark}[theorem]{Remark}

\newtheorem{remarks}[theorem]{Remarks}
\newtheorem{cor}[theorem]{Corollary}
\newtheorem{proposition}[theorem]{Proposition}

\newcommand{\bd}{\begin{definition}}
 \newcommand{\ed}{\end{definition}}

\newcommand{\bt}{\begin{theorem}}
 \newcommand{\et}{\end{theorem}}
\newcommand{\bqt}{\begin{qtheorem}}
 \newcommand{\eqt}{\end{qtheorem}}

\newcommand{\bp}{\begin{proposition}}
 \newcommand{\ep}{\end{proposition}}

\newcommand{\bl}{\begin{lemma}}
 \newcommand{\el}{\end{lemma}}
\newcommand{\bc}{\begin{cor}}
 \newcommand{\ec}{\end{cor}}

\newcommand{\bex}{\begin{example}}
 \newcommand{\eex}{\end{example}}
 
\newcommand{\bexs}{\begin{examples}}
 \newcommand{\eexs}{\end{examples}}

\newcommand{\bexe}{\begin{exercice}}
 \newcommand{\eexe}{\end{exercice}}

\newcommand{\br}{\begin{remark}}
 \newcommand{\er}{\end{remark}}
 
\newcommand{\brs}{\begin{remarks}}
 \newcommand{\ers}{\end{remarks}}

\newcommand{\bce}{\begin{center}}
\newcommand{\ece}{\end{center}}

\begin{document}
\begin{center}

{\huge On single-frequency asymptotics for  
\bigskip

the Maxwell--Bloch equations: pure states}
\bigskip\smallskip

 {\large A.I. Komech$^1$ and E.A. Kopylova}\footnote{ 
 Supported  by Austrian Science Fund (FWF)  10.55776/PAT 3476224}
 \\
{\it
Institute of Mathematics of
BOKU
University, 
\\
Gregor Mendel Strasse 33,
A-1180,
Vienna, Austria\\
}
 alexander.komech@boku.ac.at,  elena.kopylova@boku.ac.at

\end{center}

\setcounter{page}{1}
\thispagestyle{empty}
\begin{abstract}
We consider damped driven
Maxwell--Bloch equations  for  a single-mode Maxwell field
coupled to a two-level molecule. The equations 
are used for semiclassical 
description of the laser action.

Our main result is the construction of solutions with single-frequency asymptotics of the 
Maxwell field in the case of   quasiperiodic pumping.
The asymptotics hold for solutions with {\it harmonic} initial values
which are stationary states of averaged reduced equations in the interaction picture.

We calculate all harmonic states
and analyse   their stability. 
Our calculations rely on  the Hopf reduction 
 by the gauge symmetry group $U(1)$.
 The asymptotics follow by an extension of
 the  averaging  theory of Bogolyubov--Eckhaus--Sanchez-Palencia
 onto dynamical systems on manifolds. 
 The key role 
 in the application of the averaging theory
 is played by  a special a priori estimate.

  \end{abstract}
  
  \smallskip
  

\tableofcontents

\setcounter{equation}{0}
\section{Introduction}
The Maxwell--Bloch equations (MBE) were introduced by Lamb \ci{L1964} for  semiclassical 
description of the laser action \ci{AE1987,BM2011,H1984,SSL1978,SZ1997,S1986,S2012}.
The equations
are obtained as the Galerkin approximation of the Maxwell--Schr\"odinger  system
\ci{BT2009,GNS1995,K2019phys,K2013,K2022,KK2020,NW2007,S2006}.
 Our main goal is construction of solutions with single-frequency asymptotics of the Maxwell amplitude. 
 The asymptotics provisionally correspond to 
 the laser coherent radiation, which
remains a key mystery of laser action since its discovery around 1960.
The damped driven  MBE for pure states read 
\be\la{HMB2} 
\left\{\ba{rcl}
\dot A(t)\!\!&\!\!=\!\!&\!\!B(t),\quad
\dot B(t)=-\Om^2 A(t)\!-\!\ga  B(t)+cj(t)
\\\\
 i\hbar\dot C_{1}(t)\!\!&\!\!=\!\!&\!\!\hbar\om_1C_{1}(t)+ia(t)\,
  C_{2}(t),\quad
i\hbar\dot C_{2}(t)=\hbar \om_2C_{2}(t)-ia(t)\,
  C_{1}(t)
 \ea\right|,\qquad t\ge 0.
\ee
Here $A(t),B(t)\in\R$, while $C_1(t),C_2(t)\in\C$;
$\Om>0$ is the resonance frequency, $\ga>0$ is the dissipation coefficient, 
$c$~is the speed of light,  $\hbar$ is the Planck constant,
and 
 $\hbar\om_2>\hbar\om_1$ are the energy levels of active molecules.
 The current $j(t)$ and the function $a(t)$ are given by
\be\la{JOm}
j(t)=
{2}\vka
\,
 \rIm[\ov C_{1}(t)C_{2}(t)],
 \quad
a(t):=\fr {\vka}c [A(t)+A^e(t)],\qquad \vka=p\,{{\om}},\quad {\om}=\om_2-\om_1>0,
\ee
where  $p\in\R$ is proportional to the molecular dipole moment.
We 
suppose that the external field $A^e(t)$ (the pumping)   is a quasiperiodic function:
\be\la{qpp}
A^e(t)=
\rRe[ \bA^e e^{-i\Om t}+\sum_1^N \bA^e_k e^{-i\Om_k t}],
\qquad{\rm where}\qquad \bA^e,\bA^e_k\in\C,\quad \Om_k\in\R\quad{\rm and}\quad\Om_k\ne\Om.
\ee
We use the Heaviside--Lorentz units, and  in Appendix \ref{sMB} we
comment on the introduction of the equations.

For the MBE,
the charge conservation holds  $|C_{1}(t)|^2+|C_{2}(t)|^2=\const$.  We will consider solutions with $\const=1$:
\be\la{Blcc}
|C_{1}(t)|^2+|C_{2}(t)|^2=1,\qquad t\ge 0.
\ee
 Accordingly, the phase space of the system  is
$
\aX=\R^2\times \rS^3,
$
where $\rS^3$ is the unit sphere in $\C^2$.  For 
any $r>0$, 
solutions $X(t)=(A(t),B(t),C_1(t),C_2(t))$ to  (\ref{HMB2}) with $|p|/\ga=r$
 satisfy the  following a priori bounds: 
   \be\la{apri}
|X(t)|\le D_r(|X(0)|),\qquad t\ge 0,
\ee
which is proved in  Appendix \ref{swp}.
The bound  implies
  the well-posedness of the MBE
  in the phase space $\aX$.
  \smallskip

Our main goal is  asymptotics of the Maxwell amplitudes $A(t)$ and $B(t)$
    for solutions $X(t)$ to the  MBE
  as
  $p,\ga\to 0$.
  Define   complex Maxwell amplitudes by $M(t)=A(t)+iB(t)/\Om$.
The first two equations of (\ref{HMB2})  are equivalent to
\be\la{Max}
\dot M(t)=-i \Om M(t)-i\ga  M_2(t)+2icp\om\,\rIm[\ov C_1(t)C_2(t)]/\Om,\quad{\rm where}\quad 
 M_2(t)=\rIm M(t).
\ee
Note that the parameters $|p|,\ga$ are very small for many types of lasers, see Appendix \ref{aD}. 
For $p=\ga=0$,  all solutions to (\ref{Max})
are single-frequency: $M(t)=e^{-i\Om t}M(0)$.
For small $|p|$ and $\ga$, equation (\ref{Max}) implies that 
$
M(t)= e^{-i\Om t}M(0) +\ds\int_0^t R(s)ds,
$
 where $\sup\limits_{s\ge 0} |R(s)|=\cO(|p|+\ga )$  by  (\ref{apri}). Hence,  $M(t)= e^{-i\Om t}M(0) +\cO(|p|t)$
for $p/\ga=r$ with arbitrary  fixed $r\ne 0$.
In particular, for solutions with any fixed initial states $M(0)$,
 \be\la{intriad0}
 \max_{t\in [0,|p|^{-1/2}]}|M(t)-e^{-i\Om t}M(0)|=\cO({|p|^{1/2}}),
 \qquad { p\to 0},\quad{p}/\ga=r.\quad\qquad\qquad\qquad\qquad
   \ee 
   Our main results show that in the resonance case, when
  $\Om=\om$,
the time scale $|p|^{-1/2}$ in this asymptotics can be extended 
to $|p|^{-1}$ and even more in the case of special  {\it harmonic states}
of the dynamics (\ref{HMB2}).

  The  harmonic states are defined via
 the reduction of  dynamics (\ref{HMB2}) by the symmetry gauge group $U(1)$ with the action 
 \be\la{U1}
g(\theta)(A,B,C_1,C_2)=(A,B,e^{i\theta}C_1,e^{i\theta}C_2),\qquad \theta\in[0,2\pi].
\ee
The action turns the  sphere $\rS^3$, corresponding to
(\ref{Blcc}), into the Hopf fibration with the base $\rS^2$.
We introduce suitable ``global" coordinate $Q\in\C$ in the base,
 using
 the stereographic projection $\rS^2\to \C$:
   \be\la{Spi}
 Q=\fr Z{\fr 12-\bZ_3},\qquad 
\bZ_3=\pm\sqrt{\fr14-|Z|^2},\qquad Z=\ov C_1 C_2,
 \ee
 {
 where $|Z|\le 1/2$ by (\ref{Blcc}), and
 $Z=|C_1C_2|e^{i\de}$ with $\de=\theta_2-\theta_1$, where $\theta_k=\arg C_k$. 
 The base $\rS^2$ is identified with the graph of the two-valued function $Q(Z)$ on the disk $\bD:=\{Z\in\C:|Z|\le 1/2\}$.
 The coordinates $Z$ and $Q$ are invariant with respect to the action (\ref{U1}).
 The population inversion $I:=|C_2|^2-|C_1|^2$ admits the representation
  \be\la{ISi}
  I=\fr{|Q|^2-1}{|Q|^2+1},\qquad Q\in\C.
  \ee

 }
Thus, the reduction map $\Pi$
transforms  the system (\ref{HMB2}) with the phase space $\aX=\R^2\times S^3$, into the system of 
ODEs  with the factor-phase space  $\aY=\Pi \aX=\R^2\times \rS^2$.
 Accordingly, the reduced dynamics becomes the system of ODEs for the complex Maxwell amplitude
 $M(t)\in\C$
 and the complex variable $Q(t)\in\C$ 
which represents the amplitudes $(C_1(t),C_2(t))\in \rS^3$.
   The reduced 
 system can be written as
 \be\la{HMB2red}
\dot Y(t)=F(Y(t),t),\qquad t\ge 0;\qquad Y(t)=(M(t),Q(t)).
\ee
It is crucially important that the Maxwell amplitudes $M(t)$ are not changed 
by the reduction, so their properties for the Maxwell--Bloch dynamics
(\ref{HMB2}) and for the reduced  dynamics (\ref{HMB2red}) are identical.
  We consider the 
 {\it  interaction picture} of the reduced dynamics (\ref{HMB2red})
 and 
 the corresponding {\it averaged equations} with the structure
 \be\la{aver}
 \dot \bY(t)=p\bF_r(\bY(t)),\qquad t\ge 0,
 \qquad r=p/\ga.
  \ee
  We define 
  the harmonic states $\bY\in\aY$ of the reduced dynamics 
  (\ref{HMB2red}) as stationary solutions of (\ref{aver}).
   We call  $X\in\aX$ as harmonic states of the  MBE if 
   $\bY=\Pi X$ is the  harmonic state of (\ref{HMB2red}).  
We   show that  the  harmonic states depend only on the quotient $p/\ga$.
\smallskip

Our main results are  asymptotics for solutions $X(t)$ to equations (\ref{HMB2}) with a
fixed quotient $p/\ga=r$. Everywhere below we  consider the case $p>0$ only  since the extension
to $p<0$ is obvious.
{

\bt\la{t1} { Let $\Om=\om$ and the pumping be quasiperiodic. 
Then for any $r> 0$,
the following asymptotics hold.}
\smallskip\\
i)  Let $ X(0)=X^r$ be  a harmonic state of  (\ref{HMB2}), and
$\bY^r=(\bM^r,\bQ^r):=\Pi X^r$ be  the corresponding stationary solution to (\ref{aver}).
Then for solutions $(M(t),Q(t))=\Pi X(t)$ to the reduced system (\ref{HMB2red}),
 the {\it adiabatic asymptotics} holds:
 \be\la{intriad}
 \max_{t\in [0,p^{-1}]}\Big[|M(t)-e^{-i\Om t}\bM^r|{+|Q(t)-e^{-i\Om t}\bQ^r|}\Big]=\cO({p^{1/2}}),
 \qquad { p\to 0}.\quad\qquad\qquad\qquad\qquad
   \ee 
 ii)   Let, additionally,  
  $\bY^r=(\bM^r,\bQ^r)$ 
 be an
  {\it asymptotically  stable} stationary  solution to (\ref{aver}),
  and $D\subset\C^2$ be a bounded domain of attraction.
 Then for initial states $X(0)$ with $\Pi X(0)\in D$,
    \be\la{intriads}
 \max_{t\in [0,p^{-1}]}\Big[|M(t)-e^{-i\Om t}\bM(t)|{+|Q(t)-e^{-i\Om t}\bQ(t)|}\Big]=\cO({p^{1/2}}),
 \qquad { p\to 0},\qquad\qquad\qquad\qquad
   \ee 
   where $\bM(t)$ and $\bQ(t)$ satisfy 
   \be\la{bMM}
  ( \bM(t),\bA(t))\to (\bM^r,\bQ^r), \quad t\to\infty; \qquad \sup_{t\ge 0}[|\dot\bM(t)|+|\dot\bQ(t)|]=\cO( p),\qquad { p\to 0}.
   \ee
   iii) Let, additionally,  $\bY^r=(\bM^r,\bQ^r)$
 be a {\it linearly stable} 
stationary  solution to (\ref{aver}).
 Then   for initial states $X(0)=X^r$ with $\Pi X^r=\bY^r$,
  the  {\it uniform in time} asymptotics holds,
 \be\la{intri}
 \quad
 \sup_{t\in [0,\infty)}\Big[|M(t)-e^{-i\Om t}\bM^r|{+|Q(t)-e^{-i\Om t}\bQ^r|}\Big]=\cO(p),
 \qquad { p\to 0}.
   \ee 
Moreover,  for initial states $X(0)$ sufficiently close to $X^{r}$,
the  following asymptotics with attraction holds: for sufficiently
          small $\ve>0$, and some $\mu>0$,
\be\la{intriat}
 |M(t)-e^{-i\Om t}\bM^r|{+|Q(t)-e^{-i\Om t}\bQ^r|}\le C[p+d_0 e^{-p\mu t}],\qquad t\ge 0,
 \qquad{ p\le \ve},
   \ee 
   where $d_0= |X(0)-X^{r}|$.
      \et
      }
 
  \br\rm
By (\ref{ISi}), the asymptotics of type  (\ref{intriad})--(\ref{intriat}) hold also for $I(t)-I^r$,
where \linebreak $I^r=(\bQ^r|^2-1)/(|\bQ^r|^2+1)$.    
  \er
   \br\rm
   Without 
   the constraint $p/\ga=r$, the
       a priori bounds (\ref{apri}) and the asymptotics (\ref{intriad})--(\ref{intriat}) do not hold.
  For example,  the bounds must be different for small and big $|p|$ for a given dissipation $\ga>0$.
  The asymptotics 
   do not hold without this restriction since the limiting amplitudes $\bM^r$ depend on $r$.

   \er

Let us comment on
 our approach. 
  To prove the asymptotics (\ref{intriad}), we show
 that the reduced system (\ref{HMB2red}) 
 admits solutions
 $M(t)=e^{-i\Om t}\aM^r(t)$, $Q(t)=e^{-i\om t}\aQ^r(t)$
with {\it slowly varying} enveloping amplitudes $\aM^r(t)$ and $\aQ^r(t)$
for small $p$ and $\ga$ with $p/\ga=r$.
 The amplitudes are solutions to  the corresponding dynamical system
 which is  the {\it interaction picture} (or ``rotating frame representation") of (\ref{HMB2}).
  The slow variation of the amplitudes
is equivalent to the fact 
 that
the initial state $(\aM(0),\aQ(0))=(M(0),Q(0))$
is a harmonic state, i.e., a stationary state of (\ref{aver}).

 We
  calculate all  harmonic states $\bY^r=(\bM^r,\bQ^r)\in\aY$ 
  which are stationary solutions to (\ref{aver}).
 We show that 
   the states with $\bM^r\ne 0$ exist only 
  in the resonance case $\Om=\om$ and only if $\bA^e\ne 0$.
  In particular, for  single-frequency pumping
 $A^e(t)=\rRe [\bA_p e^{-i\om_p t}]$,
 the asymptotics  (\ref{intriad})--(\ref{intriat}) with $\bM^r\ne 0$ hold only in the case 
 $\bA_p\ne 0$ and 
 triple resonance
 \be\la{trip} 
  \Om=\om=\om_p.
  \ee
  We linearise the  averaged dynamics (\ref{aver}) at the harmonic states and 
  calculate the spectrum of the linearisations.
  In particular, we show that in the case  $cr>|\bA^e|$,
there exist
{\it linearly stable} 
 harmonic  states $\bY^r=(\bM^r,\bQ^r)$ of (\ref{aver}). 
 
For the proof of
  the asymptotics (\ref{intriad})--(\ref{intriat}), we extend 
 the  averaging theory of Bogolyubov type \ci[Theorem 4.3.6]{SVM2007}
 and the averaging  theory of  Eckhaus and Sanchez-Palencia 
 \ci[Theorem 5.5.1]{SVM2007} onto dynamical systems on manifolds. 
 The extension is done in Appendix \ref{aC}. The key role 
 in the application of the averaging theory
 is played by  the special a priori estimate (\ref{apri}).

   \smallskip

\br\rm
i)
 The map $\Pi:\aX\to\aY$ reduces one degree of freedom in the MBE.
 However, the role of the reduction is more significant. Namely, 
 the 
 single-frequency asymptotics of type (\ref{intriad})--(\ref{intriat})  apparently {\it do not hold}
  for  $C_1(t)$ and $C_2(t)$. 
   This conjecture is suggested by the Rabi perturbative solution \ci[(5.2.14)--(5.2.16)]{SZ1997}
      since it can
 contain incommensurable  frequencies.
\smallskip\\
ii)  Our approach relies on the averaging theory which
neglects oscillating terms. So, it gives a justification of the ``rotating wave approximation",  which is widely used in  Quantum
Optics \ci{AE1987,H1984,SSL1978,SZ1997,S1986,S2012}.
The asymptotics (\ref{intriad})--(\ref{intriat}) 
specify the  time scale and  the  error of such approximations.

\er

Let us comment on related results.
The problem of existence of time-periodic solutions 
to the MBE
has been
 discussed since 1960s. 
 The first results in this direction were obtained recently in
  \ci{CLV2019} and \ci{WWG2018}
 for various versions
 of the MBE. 
   In  \ci{WWG2018},
  the  N-th order time-periodic
  solutions 
  were constructed by perturbation techniques.
    For the phenomenological model \ci{A1985,AB1965},
  time-periodic solutions
  were constructed in \ci{CLV2019} in the absence of
   time-periodic pumping  for small interaction constants.
   The solutions are obtained as the result of a bifurcation
   relying on  homotopy invariance of the degree  \ci{CLN2017} and developing the 
   averaging arguments  \ci{CL2018}. 
     The period is determined by the bifurcation.

In \ci{K2024}, we have established the existence of 
solutions with $T$-periodic Maxwell amplitude
for any $T$-periodic pumping 
without smallness conditions.

Up to our knowledge,
the single-frequency asymptotics  for the MBE were
not constructed till now.

\smallskip

Let us comment on our exposition.
In Sections \ref{sZ} and \ref{sS} we construct an appropriate model of the reduced dynamics.
The averaged equations in the interaction picture are calculated 
 in Section \ref{sav},
and  all stationary states of these equations are calculated in
Section \ref{sst}.
In Section \ref{sstab} we analyse the stability of the stationary states. 
In Section \ref{ssin} we prove the single-frequency asymptotics 
(\ref{intriad})--(\ref{intriat}).
In Appendix \ref{swp} we establish the bound (\ref{apri}).
We comment on the introduction  of the MBE in
Appendix \ref{sMB}.
 In Appendix \ref{aD} we discuss possible treatment of the laser threshold
 and laser amplification relying on our results.

 {\bf Acknowledgements.} The authors thank S. Kuksin, M.I. Petelin, A. Shnirelman and H. Spohn
 for longterm fruitful discussions, the reviewer of our paper for useful remarks, 
 and the  Institute of Mathematics of BOKU University
for the support and hospitality. The research is supported by
Austrian Science Fund (FWF) 10.55776/PAT3476224.

\setcounter{equation}{0}
\section{$U(1)$-symmetry and   the Hopf fibration}\la{sZ}
Recall that 
the phase space of the MBE
is
$\aX=\R^2\times \rS^3$ due to
the charge conservation (\ref{Blcc}).
The system 
is
 $U(1)$-invariant with respect
to the action (\ref{U1}).
This means that the function $g(\theta)X(t)$ is a solution if $X(t)$ is. This is obvious
from
(\ref{HMB2}) and (\ref{JOm}).
  Let us denote the factor space $\aY=\aX/U(1)$
  which is the space of all orbits $X_*=\{g(\theta)X:\theta\in[0,2\pi]\}$ with $X\in\aX$.
  Thus, the reduction map  $\Pi:\aX\to\aY$ is defined by
  \be\la{Hr}
  \Pi:X=(A,B,C)\mapsto X_*=(A,B,C_*)\in\aY,\qquad X\in\aX,
  \ee
  where $C_*=P C$, and
the map $P:C\mapsto C_*$ is the Hopf fibration $\rS^3\to \rS^2$.
So, $\aY=\R^2\times \rS^2$.
The MBE induce the corresponding reduced dynamics in the factorspace
 $\aY$ which can be written as (\ref{HMB2red}).

\bl\la{lwpr}
   The reduced dynamics (\ref{HMB2red}) admits a unique global solution
   $Y(t)=(A(t),B(t), C_*(t))$ for every initial state $Y(0)\in\aY$.
\el
\bpr
The uniqueness follows from the smoothness and continuity of $F(\cdot, t)$.
To prove the existence, take any point $X_0\in \Pi^{-1}Y(0)$ and set
$Y(t)=\Pi  X(t)$, where $X(t)$ is the solution to the MBE  with the initial state $X_0$. 
Then $Y(t)$ is the solution to (\ref{HMB2red}) by definition of the reduced dynamics,
and the initial state $Y(0)=\Pi X(0)=Y_0$.
\epr
\noindent{\bf The model of the Hopf projection.}
  The Hopf projection 
 $P$
 can be represented by the map 
\be\la{z0} 
P:C=(C_1,C_2)\mapsto Z=\ov C_1 C_2
\ee
since it is constant on the fibers $\{e^{-i\theta}(C_1,C_2):\theta\in[0,2\pi]\}$. 
The orbit populations $|C_1|^2$ and $|C_2|^2$ can be expressed 
in $Z$. Indeed, $|C_1|^2+|C_2|^2=1$ and $|Z|=|C_1||C_2|$. Hence,
 the populations and the population inversion are given by
\be\la{popul2}
 |C_1|^2=\fr12\mp \sqrt{\fr14-|Z|^2},\quad
 |C_2|^2=\fr12\pm \sqrt{\fr14-|Z|^2},\quad I:=|C_2|^2- |C_1|^2=\pm\sqrt{1-4|Z|^2}.
 \ee
 All functions $Z$, $|C_1|$, $|C_2|$ and
 the population inversion $I=|C_2|^2- |C_1|^2$  are
constant on the fibers, so they are  functions on $\aY$.
Hence,
the function (\ref{z0})
defines the map $\rS^2\to \bD$ onto 
the disk $\bD:=\{Z\in\C:|Z|\le \fr 12\}$. So, $Z$
 is a local coordinate  on the base $\rS^2$ outside the points with $|Z|=\fr 12$. 
On the other hand, the formulas (\ref{popul2})
show that the map is not a homeomorphism $\rS^2\to \bD$ 
 since the inverse map 
$\bD\to \rS^2$ is two-valued.
 In particular, $Z=0$ for the Hopf projection  of ground state $(1,0)$  and for the 
projection  of the excited state $(0,1)$.

Let us introduce an appropriate model $\bS^2$ of the base $\rS^2$ with a two-folded projection
$\bS^2\to \bD$.

Denote by $\bB$ the ball $\{|\bZ|\le 1/2\}\subset\R^3$.
We identify the base $\rS^2$ with
 the sphere  $\bS^2=\{\bZ\in\R^3:|\bZ|=\fr 12\}$
and $\bD$ with the equatorial section $\bB\cap\{\bZ_3=0\}$.
 Denote by $\bS^2_{\pm}$  the hemispheres
 \be\la{mode} \bS^2_{\pm}=\{\bZ
 \in \bS^2: \bZ_3\gtrless 0\},
 \ee
  so
the orthogonal projections of $\bS^2_{\pm}$ onto the disc $\bD$
coincide.
Let us identify 
 $Z=\bZ_1+i\bZ_2$, see Fig. \ref{F2}.

\noindent{To fix the model,} we identify the North Pole $(0,0,1/2)$ with  $P(0,1)$
and the South Pole $(0,0,-1/2)$ with $P(1,0)$.
Equivalently,
 \be\la{popul4}
 \bS^2_{\pm}=
\{\bZ\in \bS^2: I\gtrless 0\}.
\ee
\br\la{rglo}\rm
In the next section, we will construct a ``global coordinate" on the base  $\bS^2$
 except one point.

\er

\noindent{\bf The reduced dynamics.}
Let us  construct suitable representation for the reduced dynamics 
(\ref{HMB2red}) in the coordinate $Z$.
Differentiating $Z(t)=P C(t)$, we obtain
from
the MBE   that
\be\la{ItZ}
\dot Z(t)=\dot{\ov C}_1C_2+\ov C_1\dot C_2=\ov{(-i\om_1 C_1+\fr a\hbar C_2)}C_2
+\ov C_1(-i\om_2 C_2-\fr a\hbar C_1)=-i{\om}\,\,Z+\fr a\hbar (|C_2|^2-|C_1|^2).
\ee
Using (\ref{JOm}) and (\ref{popul2}), we rewrite equation (\ref{ItZ}) as
 \be\la{ItZb}
\dot Z(t)=-i{\om}\,\,Z \pm b [A(t)+A^e(t)]\sqrt{1-4|Z|^2},\qquad b=\fr {\vka}{c\hbar}.
\ee

 \br\la{rsi}\rm
 i)  The right hand side of the resulting equation (\ref{ItZb})
   is function of $Z$  because of 
 $U(1)$-invariance of the system (\ref{HMB2}).
  \smallskip\\
 ii) By (\ref{popul4}) and (\ref{mode}), we must take the upper sign in  
  (\ref{popul2}) and
 (\ref{ItZb})
in the case $P C\in \bS^2_{+}$ and the lower sign in the case $P C\in \bS^2_{-}$.

 \er

\setcounter{equation}{0}
\section{Resolution of singularity and the reduced dynamics}\la{sS}

Note that the vector field $F(\cdot,t)$  in (\ref{HMB2red}) is smooth on $\aY$. On the other hand, 
the coefficients of the equation (\ref{ItZb})  are singular at $|Z|=\fr 12$.
The singularity  is due
 to the fact that the map (\ref{z0})
 is not transversal at the points $(C_1,C_2)\in \rS^3$ with $|C_1|^2=|C_2|^2=\fr 12$. 
So, we must construct a transversal modification of the map. We will do it
using the stereographic projection in the model $\bS^2$ constructed above, see (\ref{mode}).
\smallskip\\
\noindent{\bf The stereographic projection.}
The equation (\ref{ItZb}) holds with the corresponding sign
 on $\bS^2_{\pm}$. Let us rewrite the equation
  in the stereographic projection $S$ from the North Pole $(0,0,1/2)$  of the sphere $\bS^2$ which is
   the Hopf projection of the excited state $(0,1)$. 
    We identify the
   projection of a
   point  $\bZ\in \bS^2$  with the complex number
  \be\la{Sp}
Q=\fr Z{\fr 12-\bZ_3},\qquad 
\bZ_3=\pm\sqrt{\fr14-|Z|^2} \quad{\rm for}\quad \pm\bZ_3>0;\quad Z=\bZ_1+i\bZ_2\in\bD.
 \ee
 \begin{figure}[ht!]
\begin{center}
\!\!\!\!\!\!\!\!\!\!\!\!\!\!\!\!\!\!\!\!
\hspace{0.8cm}\includegraphics[width=17cm]{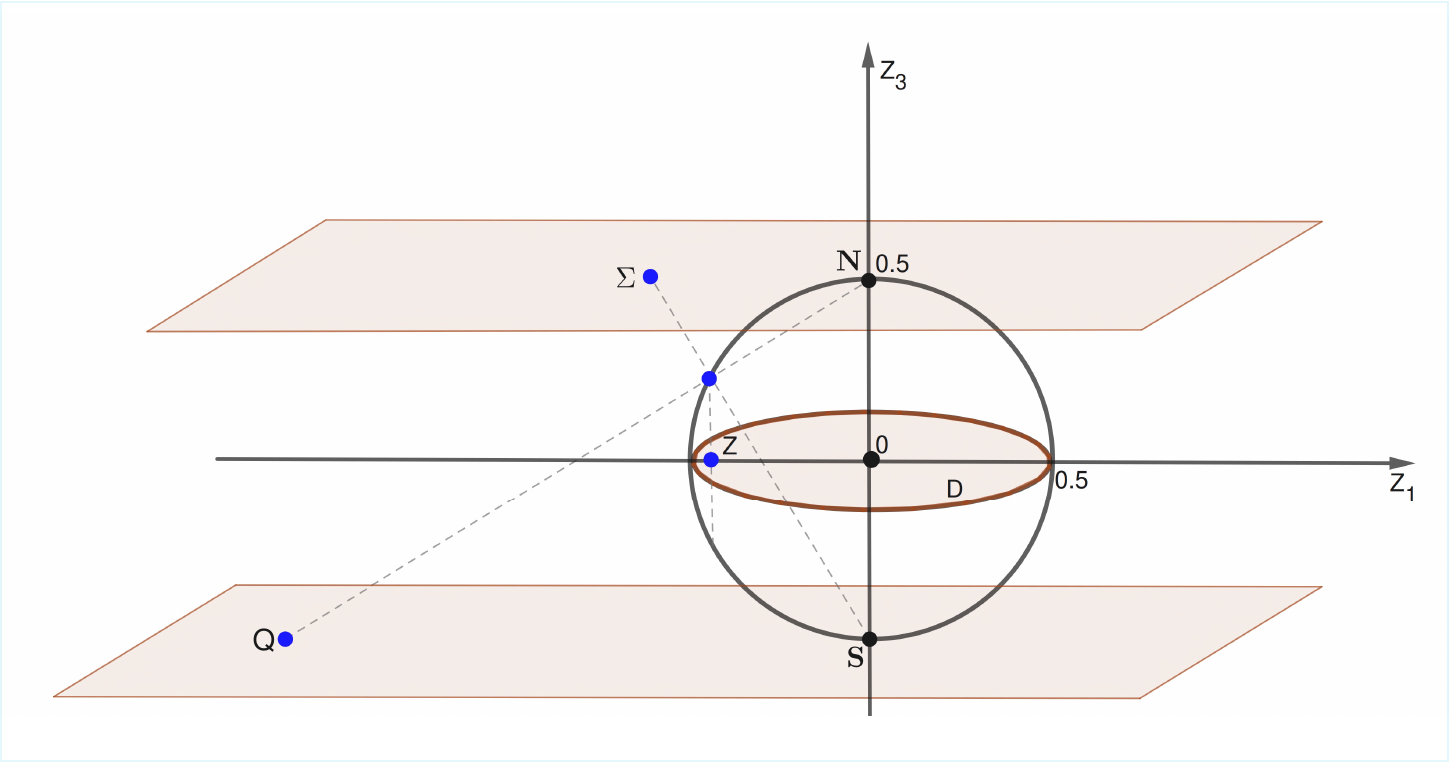}
\end{center}
\caption{
Stereografic projections}
\label{F2}
\end{figure}
    Note that  $\bZ\to Q$ is the smooth map $\dot\bS^2\to \C$,
 where $\dot\bS=\bS^2\setminus (0,0,1/2)$.
 Unlike $Z$,
 the coordinate $Q$ 
  differentiates  the hemispheres $\bS^2_{\pm}$, see Fig. \ref{F2}: by (\ref{mode}),
  \be\la{Sdiff}
  \bS^2_\pm=\{\bZ\in\bS^2: |Q|  \gtrless 1\}.
  \ee
   Hence, by (\ref{popul4}),
  \be\la{IS2}
  I\gtrless 0\Longleftrightarrow |Q|\gtrless 1.
  \ee 
 \noindent{\bf Transformation of the MBE.}
 We will express the equation (\ref{ItZb}) 
 in the variable $Q$ 
 taking into account different signs in (\ref{Sp}).
 We will see that the resulting expressions for both signs are identical.
 For any complex numbers $z_1,z_2\in\C$ we will denote
 \be\la{prod}
 z_1\cdot z_2=\rRe[\ov z_1z_2],\qquad z_1\we z_2=\rIm[\ov z_1z_2]
 \ee
 the inner and vector products of the corresponding real 2D vectors.
  Differentiating $Q(t)=Q(Z(t))$, we get from (\ref{Sp}), 
 \be\la{ItS}
 \dot Q(t)=\fr{\dot Z}{\fr12-\bZ_3}+\fr{Z}{(\fr12-\bZ_3)^2}\dot \bZ_3
 =\fr{ \dot Z}{\fr12-\bZ_3}
 \mp
 \fr{Z}{(\fr12-\bZ_3)^2}\fr{Z\cdot \dot Z}{ \sqrt{\fr14-|Z|^2}}
 =\fr{\dot Z}{\fr12-\bZ_3} -Q\fr{Q\cdot \dot Z}{\bZ_3}.
 \ee
 Therefore,  (\ref{ItZb}) can be written as
  \be\la{Itob}
 \dot Z= -i\,{\om}\,\, (\fr12-\bZ_3)Q + 2b[A(t)+A^e(t)]\bZ_3.
  \ee
Substituting this formula into the right hand side of (\ref{ItS}), we obtain
\be\la{ItSa}
 \dot Q
 =-i{\om}\,\,Q  +2b[A(t)+A^e(t)]\fr{\bZ_3}{\fr12-\bZ_3} -2 bQQ_1[A(t)+A^e(t)],\quad Q_1=\rRe Q, 
 \ee
since $Q\cdot[iQ]=0$.  Note that $\bZ_3^2+|Z|^2=\bZ_3^2+|Q|^2(\fr12-\bZ_3)^2=\fr14$ by (\ref{Sp}). Hence, 
 \be\la{Z3}
  \bZ_3=\fr12-\fr{1}{(|Q|^2+1)},\qquad
  \fr12-\bZ_3=\fr1{|Q|^2+1},\qquad
  \fr{\bZ_3}{\fr12-\bZ_3}=\fr{|Q|^2-1}2.
  \ee
 Therefore,  (\ref{ItSa}) implies that 
  \be\la{ItS1}
 \dot Q(t)=
  -i{\om} \,\,Q(t) -b[A(t)+A^e(t)](Q^2(t)+1).
  \ee 
 \br\la{rpop}\rm
 Due to  (\ref{popul2}),  the population inversion reads as (\ref{ISi}):
  \be\la{IS}
  I=\sqrt{1-4|Z|^2}=2 \bZ_3=\fr{|Q|^2-1}{|Q|^2+1}
  \ee
   since $I>0$ for $|Q|>1$ by  (\ref{IS2}).
 \er
\noindent{\bf South Pole projection.}
For any complex number $z\in\C$, we will denote $z_1=\rRe z$ and $z_2=\rIm z$.
The Maxwell amplitudes are governed by
   the first two equations of (\ref{HMB2}):
   \be\la{HMB2M}
  \dot A(t)\!=\!B(t),\,\,\,\,\,
\dot B(t)\!=\!-\!\Om^2 A(t)\!-\!\ga  B(t)+cj(t),\quad j(t)\!=\!2\vka Z_2(t)\!=\!2\vka(\fr12-\bZ_3)Q_2(t)
\!=\!2\vka\fr{Q_2(t)}{|Q(t)|^2+1}.
  \ee 
The system  (\ref{ItS1}), (\ref{HMB2M}) 
represents the smooth reduced dynamics (\ref{HMB2red}) 
in the local coordinates $(A,B,Q)$ defined in the chart $K_+=\R^2\times [\bS^2\setminus (0,0,\fr12)]$.
 Similar smooth representation   of the reduced dynamics holds in 
 the chart $K_-=\R^2\times [\bS^2\setminus (0,0,-\fr12)]$
with another
  coordinate $\Si\in\R^2$ defined as the stereographic projection from the ``South Pole" $(0,0,-\fr12)$: now (\ref{Sp}) changes to
   \be\la{Sps}
\Si(Z)=\fr Z{\fr 12+\bZ_3},\qquad 
\bZ_3=\pm\sqrt{\fr14-|Z|^2}, \quad \pm\bZ_3>0;\quad Z\in\bD.
 \ee
 Fig. \ref{F2} obviously shows
 that the North and South Pole stereografic  projections of any point of the sphere $\bS^2$ are related by
 the inversion 
\be\la{rNS}
Q=1/\ov \Si.
 \ee
 It is easy to check that
 the resulting equations 
 are very similar to (\ref{ItS1}), (\ref{HMB2M}):
     \be\la{ItS1s}
     \left\{
     \begin{split}
      \dot A(t)&=B(t),\quad
\dot B(t)\!=\!-\!\Om^2 A(t)\!-\!\ga  B(t)+cj(t),\quad j(t)=2\vka\fr{\Si_2(t)}{|\Si(t)|^2+1}
\\ 
 \dot \Si(t)&=
  i{\om}\,\, \Si(t) +b[A(t)+A^e(t)](\Si^2(t)+1)
  \end{split}\right|.
  \ee   
{\bf On representations of the reduced dynamics.}
The trajectories $Y(t)=(A(t),B(t), C_*(t))$ of the reduced dynamics
satisfy the equations  (\ref{ItS1}), (\ref{HMB2M})
in the local coordinates $(M,Q)$ on the phase space $\aY$.
The equations hold
 until $C_*(t)$ goes through the North
  Pole. Similarly, the trajectories satisfy the equations (\ref{ItS1s})
in the local coordinates $(M,\Si)$ until  $C_*(t)$ goes through the South
  Pole.

\setcounter{equation}{0}
  \section{Averaging of dynamics  in the interaction picture}\la{sav}
     Write the system (\ref{ItS1}), (\ref{HMB2M}) as
  \be\la{syst1}
  \left\{
  \ba{ll}
  \dot A( t)\!\!\!&\!\!\!=B( t),\quad
   \dot B( t)=\ds- \Om^2A( t)-\ga  B( t)+2c\vka \fr{Q_2( t)}{|Q( t)|^2+1}
   \\
\dot Q( t)\!\!\!&\!\!\!=\ds
  -i\om Q( t) - b [A( t)+A^e( t)] (Q^2( t)+1)
   \ea\right|.
  \ee 
    Let us
  denote $M( t)= A( t)+iB( t)/\Om$ and 
   rewrite  (\ref{syst1}) 
  as 
  \be\la{syst}
  \left\{\ba{rl}
  \dot M( t)&=-i\Big[\Om M( t)+\ga  M_2( t)-\ti \vka \fr{Q_2( t)}{|Q( t)|^2+1}\Big]\\\\
  \dot Q( t)&=\ds
  -i\om Q( t) -b [A( t)+A^e( t)] (Q^2( t)+1),
  \ea\right|,\qquad \ti\vka=\fr{2c\vka}\Om.
  \ee
 For small { $p,\ga>0$}, the system  is  small perturbation of the unperturbed one,
   \be\la{systu}
   \dot M( t)=-i\Om M( t),    \qquad \dot Q( t)= -i\om  Q( t).
  \ee
  {\bf The interaction picture.}
    The unperturbed system (\ref{systu}) admits 
  the single-frequency solutions
  \be\la{intr0}
   M( t)
  =e^{-i \Om t}\aM, \qquad Q( t)=e^{-i\om  t}\aQ,\qquad {\rm where}\qquad \aM,\aQ\in\C.
   \ee
Our goal is the construction of similar solutions
to the perturbed system (\ref{syst}),
 \be\la{intr}
   M( t)=e^{-i \Om t}\aM( t), \qquad Q( t)=e^{-i\om t}\aQ( t)
   \ee
with slowly varying  enveloping amplitudes: for a wide interval of time  $[0,T(p)]$,
\be\la{intr2}
  \sup_{ t\in [0,T(p)]} \Big[|\aM( t)- \aM(0)|+| \aQ( t)-\aQ(0)|\Big]\to 0,\qquad
  {  p\to 0,\quad p/\ga=r\ne 0.}
   \ee
   
   \br\la{rpi}\rm
For such solutions, $|Q( t)|\sim \const$. Therefore,  
the corresponding population inversion also is almost constant by (\ref{IS}).
 \er

Substituting (\ref{intr}) into (\ref{syst}), we obtain  the dynamical equations
for the enveloping amplitudes:

 \be\la{systim}
  \left\{\ba{rl}
  \dot \aM( t)&=-ie^{i \Om t}\Big[\ga\,\rIm( e^{-i \Om t}\aM( t)) -\ti \vka
   \fr{\rIm( e^{-i\om t}\aQ( t))}{|\aQ( t)|^2+1}\Big]\\\\
  \dot \aQ( t)&=-b e^{i\om t}[\rRe( e^{-i \Om t}\aM( t))+A^e( t)]  (\aQ^2( t)e^{-2i\om t}+1)
  \ea\right|.
  \ee
  The equations are called as the interaction picture of (\ref{syst}).
The  equations can be written as
  \be\la{perN}
  \dot \aM( t)=p\, f_r(\aM( t),\aQ( t), t),\qquad \dot \aQ( t)=p \,g_r(\aM( t),\aQ( t), t),
  \ee
 where the functions $f_r$ and $g_r$ are given by
      \be\la{fg}
  \left\{\ba{rl}
  f_r(\aM,\aQ, t)&\!\!=
 -i e^{i\Om t}
\Big[ \ga_1  (\aM_2\cos\Om  t-\aM_1\sin \Om t)-\vka_1\fr{\aQ_2\cos \om  t-\aQ_1\sin\om   t)}{|\aQ|^2+1}\Big]
\\
\\
g_r(\aM,\aQ, t)&\!\!=\ds
  -b_1e^{i\om  t}
 [\aM_1\cos  \Om t+\aM_2\sin\Om  t+A^e( t)] 
 (\aQ^2e^{-2i\om t}+1) 
\ea\right|,
  \ee
  and the parameters are
  \be\la{P1}
  \ga_1=\fr{\ga}p=\fr1{ r},\qquad \vka_1=\fr{\ti\vka}p=\fr{2c\,{\om}}{\Om},\qquad
  b_1=\fr {b}{p}=\fr{\om}{c\hbar}.
  \ee
  \br\rm
   It is important that 
  the coefficients $\ga_1,\vka_1,b_1$ depend only on 
  $\Om,\om,r$.
    Hence,
    for any fixed $\Om,\om>0$ and $r\ne 0$,
     the  asymptotics of solutions to the systems (\ref{perN})
  as { $p\to 0$ and $p/\ga=r$}
   can be calculated by methods of the averaging   theory  \ci{B1961,SVM2007}.
  \er
  \noindent{\bf The averaging.}      
  Averaging equations (\ref{perN}), we obtain the equation (\ref{aver}) in the form
    \be\la{perNa}
  \dot \bM( t)
  =p\,\ov f_r(\bM( t),\bQ( t)),\qquad \dot \bQ( t)=p\,\ov g_r(\bM( t),\bQ( t)),
  \ee
 where 
 \be\la{fga}
  \ov f_r(\bM,\bQ)=
  \langle
  f_r(\bM,\bQ,\cdot)
  \rangle=
  \lim_{T\to\infty}\fr1T \int_0^T f_r(\bM,\bQ, t)d t,\qquad
\ov g_r(\bM,\bQ)=
  \langle
  g_r(\bM,\bQ,\cdot)
  \rangle.
  \ee
  Let us calculate the averages  (\ref{fga}). The results differ drastically 
for the resonance case $\Om=\om$ and non-resonance $\Om\ne \om$.
In notation (\ref{qpp}),
    $\bA^e= \bA^e_1+i\bA^e_2$, where
\be\la{Apav}
 \bA^e
=2\langle A^e( t)e^{i  \Om t} \rangle,\qquad
  \bA^e_1
=2\langle A^e( t)\cos   \Om t \rangle,\qquad
  \bA^e_2=2\langle A^e( t)\sin  \Om   t \rangle.
  \ee

\noindent
{\bf Resonance case $\Om=\om$.} Using the expressions (\ref{fg}), we obtain
\be\la{avrnr}
\ds\ov f_r(\bM,\bQ)= -\fr{i\ga_1}2 ( \bM_2- i\bM_1)+\fr{i\vka_1}{2(|\bQ|^2+1)}(\bQ_2- i\bQ_1)
= -\fr{\ga_1}2 \bM+\fr{\vka_1}2 \fr\bQ{|\bQ|^2+1}.
\ee
Similarly, 
\beqn
\la{avr}
\nonumber
\ov g_r(\bM,\bQ)&=&-\fr{b_1}2\Big[(\bM_1-i\bM_2+ \bA^e_1-i\bA^e_2)(\bQ_1^2-\bQ_2^2+2i\bQ_1\bQ_2)+\bM+\bA^e\Big] \\
\nonumber 
&=&\fr{b_1}2\Big[(\bM_1+ \bA^e_1)(|\bQ|^2-2\bQ_1(\bQ_1+i\bQ_2))+i(\bM_2+\bA^e_2)(|\bQ|^2+2\bQ_2(i\bQ_1-\bQ_2))-\bM-\bA^e\Big]\\
\nonumber
&=&\fr{b_1}2\Big[(\bM+ \bA^e)(|\bQ|^2-1)-2(\bM_1+ \bA^e_1)\bQ_1(\bQ_1+i\bQ_2)-2(\bM_2+\bA^e_2)\bQ_2(\bQ_1+i\bQ_2)]
\nonumber\\
& =&\ds \fr {b_1}2\Big[
\bM(|\bQ|^2-1) - \ds
2 \bQ\, \,\bM\!\cdot\!\bQ
-
2\ds
   \bQ\,\,\bA^e\!\cdot\!\bQ+{ \bA^e}(|\bQ|^2-1)\Big]
 \nonumber\\
& =& \ds \fr {b_1}2\Big[
[\bM+ \bA^e](|\bQ|^2-1) - \ds
2 \bQ\, \,[\bM\!+ \bA^e]\cdot\!\bQ\Big].
   \eeqn

\noindent      
{\bf Non-resonance case $\Om\ne \om$.} In this case the calculations (\ref{avrnr}) simplify to
\be\la{avrnr2}
\ov f_r(\bM,\bQ)=
 \ds-\fr{\ga_1}2 \bM,
 \ee
which implies 
nonexistence of stationary states with nonzero Maxwell field 
for the system (\ref{perNa}),
and also
the decay $\bM( t)=\ds\bM(0)e^{-\fr{\ga} 2 t}$.
  In contrast,
the resonance equation (\ref{avrnr})  includes the interaction term 
which can prevent the decay of
 the Maxwell amplitude, that is expected physically.
 
 \noindent{\bf The averaged system.}
By (\ref{avrnr}) and (\ref{avr}), the  averaged equations (\ref{perNa})
(or (\ref{aver})) 
in the resonance case
read as
\be\la{aveq}
\!\!\!\!\!\!\left\{\ba{rcl}
\dot\bM( t)&=&
p
\Big[-\fr{\ga_1}2 \bM+\fr{\vka_1}2 \fr\bQ{|\bQ|^2+1}\Big]
\\\\
\dot\bQ( t)&=&

\ds \fr {pb_1}2\Big[
[\bM+ \bA^e](|\bQ|^2-1) - \ds
2 \bQ\, \,[\bM\!+ \bA^e]\cdot\!\bQ\Big]
\ea\right|.
\ee

      \setcounter{equation}{0}
\section{Stationary states for  averaged dynamics (harmonic states)}\la{sst}

In this section, we 
calculate all {\it harmonic states} $(\bM,\bQ)$,
i.e., stationary states 
for the averaged reduced equations (\ref{aveq}),
 in the resonance case $\Om=\om$. The states satisfy the system
\be\la{perNa1}
0=
 \ds-{\ga_1} \bM+{\vka_1} \fr\bQ{|\bQ|^2+1},
\qquad
0=
 [\bM+ \bA^e](|\bQ|^2-1) - \ds
2 \bQ\, \,[\bM\!+ \bA^e]\cdot\!\bQ.
\ee
It is important that the stationary equations depend on { $r$} but {\it do not depend} on { $p$}.
The first equation of (\ref{perNa1}) gives
\be\la{al}
\bM=\al \bQ\quad{\rm where}\quad \al=\fr{\vka_1}{\ga_1(|\bQ|^2+1)}
{ =\fr{2cr}{(|\bQ|^2+1)}>0}
\ee
according to (\ref{P1}) with $\Om={\om}$.
Now the second equation of (\ref{perNa1}) becomes
\be\la{sec}
-\ds
\al\bQ(|\bQ|^2+1) 
-
2\ds
   \bQ\, \bA^e\!\!\cdot\!\bQ+{ \bA^e}(|\bQ|^2-1) =0.
   \ee
   Note that
   in the case $\bA^e=0$, we have $\bQ=\bM=0$, so we consider below   $\bA^e\ne 0$ only.
   Taking the vector product of (\ref{sec}) with $\bQ$, we get 
$$
 \bA^e\we\bQ(|\bQ|^2-1)=0.
 $$
Hence, it suffices to consider  the cases 1) $|\bQ|=1$, and 2) $|\bQ|\ne 1$ while  $\bA^e\we\bQ=0$.
\subsection{Zero population inversion}
In the case 
 $|\bQ|=1$, we have $|Z|=1/2$, so 
by (\ref{popul2}),
the population inversion vanishes: $I=0$.
\bl\la{lS-1} 
Let $\bA^e=e^{i\varphi}|\bA^e|\ne 0$.  Then solutions $(\bM,\bQ)$ to  the system 
(\ref{perNa1}) with $|\bQ|=1$
 exist iff  
\be\la{al0}
{ { cr}\le |\bA^e|. }
 \ee
The solutions read as
 \be\la{sol}
 \bQ_\pm=e^{i(\varphi\pm \theta)},\quad  \bM_\pm=cre^{i(\varphi\pm \theta)},
  \quad{\rm where} \quad \theta=\arccos\frac{-cr}{|\bA^e|}.
\ee
For  ${ cr} < |\bA^e|$,  there are two solutions.  For $cr= |\bA^e|$, the unique solution is given by
\be\la{S1}
\bQ=-\bA^e/|\bA^e|,\qquad  \bM=-\bA^e.
\ee
\el
\bpr
For $|\bQ|=1$, equation (\ref{sec}) simplifies to
\be\la{AeS}
\bA^e\cdot\bQ=-\ds\fr{\vka_1}{2\ga_1}=-cr
\ee
according to (\ref{P1}) with $\Om={\om}$. Hence,
\be\la{seq}
|\bA^e|\cos \widehat{(\bA^e, \bQ)}=-cr
\ee
Evidently, the equation has solutions  iff  (\ref{al0}) holds.  
Finally,   (\ref{sol})  and (\ref{S1}) follow from (\ref{seq}) and  (\ref{al}).
\epr
\subsection{Nonzero population inversion}
Here we calculate stationary states of the system (\ref{perNa1})  in the case $|\bQ|\ne 1$,  $\bA^e\we\bQ=0$. 
\begin{lemma}\la{lS-2}
Let $\bA^e\ne 0$.  Then nonzero  solutions $(\bM,\bQ)$ to (\ref{perNa1}) with $|\bQ|\ne 1$ and  $\bA^e\we\bQ=0$
exist iff 
\be\la{nzD}
{ cr> |\bA^e|. }
\ee
The solutions are given by 
\be\la{nz22}
\bQ_\pm=\Big[-{cr}\pm
  \sqrt{c^2r^2-|\bA^e|^2}\Big]\frac{\bA^e}{|\bA^e|^2},
\quad  \bM_\pm=\al_\pm\bQ_\pm,\quad {\rm where}\quad\al_\pm=\fr{2cr}{|\bQ_\pm|^2+1}.
\ee
The following relations hold:
\be\la{Spm}
{ |\bQ_+|<1\quad{\rm and}\quad |\bQ_-|>1.}
\ee
\end{lemma}
\bpr   
Assume first that $ \bA_1^e\ne 0$. Then  $ \bA^e\we\bQ= 0$ implies $\bQ_2=\frac{ \bA_2^e}{ \bA_1^e}\bQ_1$.
Hence
$$
 \bA^e\!\!\cdot\!\bQ= \bA_1^e\bQ_1+ \bA_2^p\frac{ \bA_2^e}{ \bA_1^e}\bQ_1=\bQ_1\frac{|\bA^e|^2}{ \bA_1^e},
\qquad
|\bQ|^2=\bQ_1^2\Big(1+\frac{|\bA_2^e|^2}{|\bA_1^e|^2}\Big)=\bQ_1^2\frac{|\bA^e|^2}{|\bA_1^e|^2}.
$$
Substituting into \eqref{sec}, we obtain for the first component,
$$
-\fr{\vka_1}{\ga_1}\bQ_1 -2\bQ_1^2\frac{|\bA^e|^2}{ \bA_1^e}+\bQ_1^2\frac{|\bA^e|^2}{\bA_1^e}-{ \bA_1^e}=0.
$$
Simplifying, we get 
\be\la{sec2}
\frac{\ga_1|\bA^e|^2}{ \bA_1^e}\bQ_1^2+\vka_1\bQ_1+\ga_1 \bA_1 ^e=0.
\ee
Hence we get the solutions
\be\la{nz}
\bQ_\pm=\frac{-\vka_1\pm\sqrt{\vka_1^2-4\ga_1^2|\bA^e|^2}}{2\ga_1}\frac{\bA^e}{|\bA^e|^2},
\qquad \bM_\pm=\al_\pm\bQ_\pm,\quad{\rm where}\quad 
\al_\pm=\fr{\vka_1}{\ga_1(|\bQ_\pm|^2+1)}.
\ee
The formula also holds
in the case when $ \bA_1^e=0$ and $ \bA_2^e\ne 0$.
At last, the solution can be rewritten as (\ref{nz22})
since $\vka_1/\ga_1=2cr$ by (\ref{P1}) with ${\Om}=\om$.
 Finally, (\ref{Spm})  holds since $|\bQ_+|\cdot|\bQ_-|=1$ by (\ref{nz22}).
\epr

\subsection{ Fully inverted states are  not  harmonic}

Our analysis above does not include the North Pole which corresponds to $Q=\infty$.
To extend the analysis, we can use the 
equations (\ref{ItS1s}) corresponding to the  
stereographic projection from the South Pole (\ref{Sps}).
{ Slightly modifying  calculations of Section \ref{sav}, it is easy to check
 that
for $\bA^e\ne 0$ and any $r\ne 0$,
the subset $\{\Si=0,\bM\in\C\}\subset \aY=\aX/U(1)$ does  not contain  stationary states} of the corresponding 
averaged equations in the interaction picture. { The population inversion at these states is given by 
(\ref{IS}) with $\bQ=\infty$, i.e., 
$I=1$.
So,
all fully inverted states are not harmonic.}

\subsection{Summary}
The results of this section imply the following corollary. 
  
 \bc\la{r2}
 
 i) In the non-resonance case, $\Om\ne \om$, 
 the averaged system  (\ref{perNa}) 
does not have
stationary states with nonzero Maxwell field.
\smallskip\\
ii) In the resonance case, $\Om= \om$,
the system  (\ref{perNa})  admits  two stationary states
 given by (\ref{sol}) for ${ cr}<|\bA^e|$ and by
 (\ref{nz22})  for ${ cr}> |\bA^e|$.  In the case $ cr= |\bA^e|$, the  states coincide and are 
 given by (\ref{S1}).
   \ec
\br
\la{runif}
\rm
Let $(\bM(t),\bQ(t))$ be a solution to 
the averaged
equations (\ref{aveq}). Then  $e^{i\phi}(\bM(t),\bQ(t))$ with $\phi\in\R$ is the  solution
to the same equations with $\bA^e$ replaced by $e^{i\phi}\bA^e$. Hence,
this correspondence also holds for solutions to
 stationary equations (\ref{perNa1})  that
is obvious  from the formulas (\ref{sol}) and (\ref{nz22}).

\er

  \setcounter{equation}{0}
  \section{Stability of harmonic states}\la{sstab}
    
  In this section we calculate the spectra of the linearization of the averaged 
  equations (\ref{perNa}) at the harmonic states calculated above 
  in the resonance case $\Om=\om$.
  As the result, we will obtain that the states (\ref{sol}) and (\ref{S1}) with $|\bQ|=1$ are not linearly stable, the states
  $\bQ_+$ in
   (\ref{nz22})
  is linearly stable, while $\bQ_-$ 
  is unstable .
Recall the resonance formulas (\ref{avrnr}), (\ref{avr}): omitting the index $r$, 
$$
\left\{\ba{rl}
\ov f(\bM,\bQ)=\!\!\!&\!\!\!
\ds-\fr{\ga_1}2 \bM+\fr{\vka_1}2 \fr\bQ{|\bQ|^2+1}, 
 \\
 \ov g(\bM,\bQ)
 = \!\!\!&\!\!\! \ds\frac{b_1}2\Big[\bM(|\bQ|^2-1) -2\bQ\, \,\bM\!\cdot\!\bQ-2\bQ\,\,\bA^e\!\cdot\!\bQ+\bA^e(|\bQ|^2-1)\Big].
\ea\right.
$$
Recall that $z_1$ and $z_2$ denote the real and imaginary parts for any complex number $z\in\C$. 
Differentiating, we get for $i,j=1,2$,
$$
~\!\!\!\!\!\!\!\!\!\!\!
\left\{\!\!\!\!\ba{ll}
\ds\frac{\pa \ov f_i}{\pa \bM_j}=-\frac{\ga_1}2\delta_{ij},\,\,&\,\,
\ds\frac{\pa \ov f_i}{\pa  \bQ_j}=\frac{\vka_1}2\Big(\frac{\delta_{ij}}{|\bQ|^2+1}-\frac{2 \bQ_i \bQ_j}{(|\bQ|^2+1)^2}\Big),
\\\\
\ds\frac{\pa \ov g_i}{\pa  \bM_j}\!=\!\fr {b_1}{2}\delta_{ij}(|\bQ|^2\!-\!1)-b_1  \bQ_i \bQ_j,
\,\,&\,\,
\ds\frac{\pa \ov g_i}{\pa  \bQ_j}=b_1
[\bM_i\bQ_j-\delta_{ij}\bM\cdot \bQ-\bQ_i\bM_j-\delta_{ij} \bA^e\!\cdot\! \bQ -\bQ_i\bA^e_j+\bA^e_i \bQ_j].
\ea
\right.
$$
Taking into account (\ref{al}), we obtain
the Jacobian 
    \be\la{diff2}
J=\begin{pmatrix}
-\ds\frac{\ga_1}2\delta_{ij}&
\ds\frac{\vka_1}2\Big(\frac{\delta_{ij}}{|\bQ|^2+1}-\frac{2\bQ_i\bQ_j}{(|\bQ|^2+1)^2}\Big)
\\\\
\ds\fr {b_1}{2}(|\bQ|^2-1)\delta_{ij}-b_1 \bQ_i\bQ_j
&
-b_1
[\delta_{ij}\bM\cdot \bQ+ \delta_{ij} \bA^e\!\cdot\! \bQ
+
\bQ_i \bA^e_j-\bA^e_i \bQ_j]
\end{pmatrix}
\ee
 We can choose  real coordinates  in $\C=\R^2$  such that 
$\bQ=(|\bQ|,0)$. 
Then
   the Jacobian (\ref{diff2}) becomes
\be\la{mafo}
J=\begin{pmatrix}
-\fr{\ga_1}2&0&\frac{\vka_1}2\frac{1-|\bQ|^2}{(|\bQ|^2+1)^2}&0\\
0&-\fr{\ga_1}2&0&\frac{\vka_1}2\frac{1}{|\bQ|^2+1}\\
-\fr {b_1}{2}(|\bQ|^2+1)&0&
-b_1(s+\fr{\vka_1}{\ga_1}\fr{|\bQ|^2}{|\bQ|^2+1}) &-b_1w
\\
0&\fr {b_1}{2}(|\bQ|^2-1)&
 b_1w&-b_1(s+\fr{\vka_1}{\ga_1}\fr{|\bQ|^2}{|\bQ|^2+1})
\end{pmatrix},
\ee
where 
$$
s= \bA^e\cdot\bQ,\qquad w=\bA^e\we\bQ.
$$
      \subsection{Zero population inversion.}
      
      Denote $s_\pm= \bA^e\cdot\bQ_\pm$ and  $w_\pm=\bA^e\we\bQ_\pm$.
      For the states 
       (\ref{sol})
       with $|\bQ_\pm|=1$, 
  the matrix (\ref{mafo}) becomes
  \be\la{mafo2}
J_\pm=
\left(\ba{cccc}
-\fr{\ga_1}2&0&0&0\\
0&-\fr{\ga_1}2&0&\frac{\vka_1}4\\
-b_1&0& 0  &b_1w_\pm\\\
0&0&-b_1w_\pm&0
\ea\right)
\ee
since $s_\pm+\fr{\vka_1}{\ga_1}\fr{|\bQ|^2}{|\bQ|^2+1}=\bA^e\cdot\bQ+\fr{\vka_1}{2\ga_1}=0$ by \eqref{AeS}.
The eigenvalues of $J_\pm$ equal
\be\la{laze}
\lam_{1,2}=-\ga_1/2,\qquad \lam_{3,4}= \pm ib_1w_\pm.
\ee
Note that the 
spectrum of the linearised system (\ref{perNa}) at the 
stationary states  (\ref{sol}), (\ref{S1}) 
consists of $p\lam_{1,2}$ and $p\lam_{3,4}$, where $p\in\R$.
Thus, we have proved the following lemma.
\bl\la{c0}
For   ${cr} \le  |\bA^e|$ and all $p\in\R$,  stationary states (\ref{sol})
of the system (\ref{perNa})
are not linearly stable.
\el
\subsection{Nonzero population inversion}
Here we  consider  the stationary states (\re{nz22}) corresponding to  ${|cr|}> |\bA^e|$. 
For a moment, let us write $\bQ$ instead of $\bQ_\pm$.
 By   (\ref{nz22}), 
 $\bA^e=(-|\bA^e|,0)$ in the case $\bQ=(|\bQ|,0)$.   Hence,  equation
 (\ref{sec2}) becomes
$$
-{\ga_1|\bA^e|}|\bQ|^2+\vka_1|\bQ|-\ga_1 |\bA ^e|=0.
$$
This implies  that
$\ds\fr{\vka_1}{\ga_1}\fr{|\bQ|}{|\bQ|^2+1}=|\bA^e|$, 
and  therefore
$$
s+\fr{\vka_1}{\ga_1}\fr{|\bQ|^2}{|\bQ|^2+1}=
\bA^e\!\cdot\!\bQ+\fr{\vka_1}{\ga_1}\fr{|\bQ|^2}{|\bQ|^2+1}=
-|\bA^e||\bQ|+\fr{\vka_1}{\ga_1}\fr{|\bQ|^2}{|\bQ|^2+1}
=|\bQ|\Big[\fr{\vka_1}{\ga_1}\fr{|\bQ|}{|\bQ|^2+1}-|\bA^e| \Big]=0.
$$
Now $w=\bQ\we\bA^e=0$, and the Jacobian (\ref{mafo}) becomes
 \be\la{mafo3}
J=
\left(\ba{cccc}
-\fr{\ga_1}2&0&\frac{\vka_1}2\frac{1-|\bQ|^2}{(|\bQ|^2+1)^2}&0\\
0&-\fr{\ga_1}2&0&\frac{\vka_1}2\frac{1}{|\bQ|^2+1}\\
-\fr {b_1}{2}(|\bQ|^2+1)&0&
0&0
\\
0&\fr{b_1}{2}(|\bQ|^2-1)&
 0&0
\ea\right).
\ee
 The characteristic equation  $\det(2J-\lam)=0$ reads as
  $$
 -(\ga_1+\lam)\begin{bmatrix}
 -(\ga_1+\lam)   &  0  & \frac{{\vka_1}}{|\bQ|^2+1}\\
   0   &- \lam &     0\\
    b_1(|\bQ|^2\!-\!1) &0   &-\lam
  \end{bmatrix} 
  - \frac{{\vka_1}(|\bQ|^2\!-\!1)}{(|\bQ|^2+1)^2}\begin{bmatrix}
  0& -(\ga_1+\lam) &  \frac{{\vka_1}}{|\bQ|^2+1}\\
 -b_1(|\bQ|^2\!+\!1)   &0 &0\\
  0 &  b_1(|\bQ|^2\!-\!1)&-\lam
 \end{bmatrix} =0.
  $$
  Evaluating, we obtain
  $$
  (\ga_1+\lam)^2\lam^2 -  (\ga_1+\lam)\lam \frac{b_1{\vka_1}(|\bQ|^2-1)}{|\bQ|^2+1} 
  +\frac{{\vka_1}(|\bQ|^2-1)}{(|\bQ|^2+1)^2}\Big(b_1^2{\vka_1}(|\bQ|^2-1)- b_1(|\bQ|^2+1) (\ga_1+\lam)\lam \Big)=0.
  $$
  Equivalently, 
  $$
  (\ga_1+\lam)^2\lam^2 -  2(\ga_1+\lam)\lam \frac{b_1{\vka_1}(|\bQ|^2-1)}{|\bQ|^2+1}  +\frac{b_1^2{\vka_1}^2(|\bQ|^2-1)^2}{(|\bQ|^2+1)^2}=
  \Big[
  (\ga_1+\lam)\lam-b_1{\vka_1}I
  \Big]^2
  =0
  $$
   by (\ref{IS}).   
 Hence, we have the following roots of multiplicity two:
\be\la{lam1}
\qquad\qquad\qquad\quad \quad
\lam_{1,2}=\lam_{1,2}(\bQ)=\ds\frac12\Big[{-\ga_1\pm\sqrt{\ga_1^2+4b_1{\vka_1}I(\bQ)}}\Big],
\qquad I(\bQ)=\fr{|\bQ|^2-1}{|\bQ|^2+1}.
\ee
By (\ref{P1}),  the 
spectrum of the linearised system (\ref{perNa}) at the 
stationary states  (\ref{nz22})
consists of 
$$
p\lam_{1,2}=\ds\frac12\Big[{-\ga \pm \sqrt{{\ga^2}+4p^2b_1{\vka_1}I(\bQ)}}\Big].
$$
Note that 
 $\vka_1,b_1>0$ by (\ref{P1}),
 so for $p>0$,
 \be\la{lam12}
 \rRe[\lam_{1,2}(\bQ)]=-\nu<0
 \ee
iff
$I(\bQ)<0.$
 Recall that
$I(\bQ_+)<0$ and $I(\bQ_-)>0$  by  (\ref{Spm}).
%
%
%
Thus, we have proved the following lemma.
 
  \bl\la{cex} 
  In the case $cr>|\bA^e|$,
stationary states $(\bM_+,\bQ_+)$ resp. $(\bM_-,\bQ_-)$
of the system (\ref{perNa}), which are 
given by  (\ref{nz22}),  are linearly 
stable resp. unstable.

 \el

        \setcounter{equation}{0}
  \section{Single-frequency asymptotics}\la{ssin}
  
   In this section we prove 
    asymptotics (\ref{intriad})--(\ref{intriat}) in the resonance case $\Om=\om$.
  We obtain the asymptotics
     applying the results of the  averaging theory \ci{SVM2007}. 
   To justify the application, we are going to check suitable properties of the 
   system (\ref{perN}).

   \subsection{KBM vector field}
 Let us denote $v_r(\aM,\aQ,t)=(f_r(\aM,\aQ,t), g_r(\aM,\aQ,t))$ 
 the vector field (\ref{fg})
   of the system  (\ref{perN}).
 It is easy to check that in  the case of  {\it almost periodic pumping} $A^e(t)$ the asymptotics holds
  \be\la{KBM}
  \sup_{(\bM,\bQ)\in \cD}\fr1{T}
 \Big| \int_0^T[v_r(\bM,\bQ, t)-\ov v_r(\bM,\bQ)]d t\Big|
  \to 0,\qquad T\to \infty
  \ee  
 for each
 $r\!\ne\!0$
  and any bounded   region $\cD\subset\C^2$.
 Hence, according to \ci[Definition 4.2.4]{SVM2007}, $v_r$
   is a
KBM (Krylov--Bogolyubov--Mitropolsky) vector field on any bounded   region $\cD\subset\C^2$.
Furthermore,    for the {\it quasiperiodic pumping} (\ref{qpp}),
formulas
 (\ref{fg}) imply that
    \be\la{KBM2}
  \de_\cD(p):=p
  \sup_{(\bM,\bQ)\in \cD}
  \,\,
  \sup_{T\in[0,p^{-1}]}
 \Big| \int_0^T[v_r(\bM,\bQ, t)\!-\!\ov v_r(\bM,\bQ)]d t\Big|
=\cO(p),\,\,\, p\to 0,
  \ee    
  where $\de_\cD(p)$
  is  
the corresponding {\it order function} defined in \ci[Lemma 4.6.4]{SVM2007}.

  \br\la{rD}\rm
     For more general almost periodic pumping $A^e(t)$, the order function
  can be different \ci[Section 4.6]{SVM2007}.
  \er
   \subsection{Single-frequency  asymptotics in the interaction picture}\la{sing}

Here we prove Theorem \ref{t1}.    We  denote by  $\aY(t)=(\aM(t),\aQ(t))$ and $\bY(t)=(\bM(t),\bQ(t))$  the solutions to the 
 systems (\ref{perN}) and (\ref{aveq}), respectively, with the same initial value $\aY(0)=\bY(0)=\Pi X(0)$.
\smallskip\\
{
i)  Let us prove the asymptotics 
 \be\la{intriadeq}
 \max_{t\in [0,p^{-1}]}\Big[|\aM(t)-\bM^r|+|\aQ(t)-\bQ^r|\Big]=\cO({p^{1/2}}),
 \qquad { p\to 0},\quad{p}/\ga=r.\quad\qquad\qquad\qquad\qquad
   \ee 
   In other words,
      \be\la{Bt521}
 \quad\,\,\,\,\,\max_{ t\in[0,p^{-1}]} |  \aY( t)-\bY^r|= \cO({p^{1/2}}),
 \qquad p\to 0,\quad{p}/\ga=r,
  \ee 
   where $\bY^r=(\bM^r,\bQ^r)$ is a stationary state for the averaged system  (\ref{aver}), or equivalently,  to (\ref{aveq}).
Note that  
   $\bY^r:=\Pi X^r=(M(0),Q(0))$  since $X(0)=X^r$. Hence, (\ref{intr}) implies that
  $\aY(0)=(M(0),Q(0))= \bY^r$. Recall that $v_r$
  is the KBM vector field  in each bounded region $\cD\subset\C$.
    Hence, (\ref{Bt521}) would follow from Theorem 4.3.6  of \ci{SVM2007} for the case when
    all trajectories $\aY(t)$ are uniformly bounded for small $p>0$.

   The Maxwell amplitudes $\aM(t)$ are uniformly bounded
  by the a priori estimate (\ref{apr2}) which is uniform in small $p>0$ and $\ga=pr$
  for every fixed $r>0$. On the other hand, the trajectories $\aQ(t)$ can be unbounded
  if  $\bZ(t)$, 
  defined by (\ref{Sp}),
  approaches the North pole $(0,0,1/2)$ of the sphere $\rS^2$. 
  This is why the  Theorem 4.3.6  of \ci{SVM2007} cannot be applied since the vector field
   $v_r(\aM,\aQ)$ is not  Lipschitz continuous in $\aQ\in\C$.

  This suggests   to consider  (\ref{perN}) 
  as the dynamical system 
  for  $(\aM(t),\bZ(t))\in\C\times\rS^2$ 
   with
 the corresponding vector field $\bv_r(\aM,\bZ)$ on the manifold $\C\times\rS^2$.
The key observation is that the vector field $\bv_r(\aM,\bZ)$
  is Lipschitz continuous  on every compact subset  of the manifold
  $ \C\times\rS^2$ in the sense (\ref{LIP}).
    This follows from the representations (\ref{fg}) in the coordinates $\aM,\aQ$
  which holds outside  the North pole,
  and similar representations for the system corresponding to (\ref{ItS1s})
  outside  the South pole.

Hence,
  the bounds (\ref{Bt521}) hold by Theorem \ref{tC} which 
   extend 
  the Theorem 4.3.6  of \ci{SVM2007} onto  dynamical systems on compact  manifolds.
  As the result, (\ref{intriadeq}) follows, so (\ref{intriad}) holds by
the relations (\ref{intr}).
  \smallskip\\
    ii) By the same Theorem \ref{tC},
 \be\la{Bt522}
 \quad\,\,\,\,\,\max_{ t\in[0,p^{-1}]} | \aY(t)-\bY(t)|=
 \cO({p^{1/2}}),
 \qquad p\to 0,\quad{p}/\ga=r.
  \ee 
In other words,
  \be\la{intriadseq}
 \max_{t\in [0,p^{-1}]}\Big[|\aM(t)-\bM(t)|+|\aQ(t)-\bQ(t)|\Big]=\cO({p^{1/2}}),
 \qquad { p\to 0},\quad{p}/\ga=r. \qquad\qquad\qquad\qquad
   \ee 
Hence, asymptotics (\ref{intriads}) follows by  (\ref{intr}).
Finally, (\ref{bMM}) holds since $\Pi X(0)\in D$ and $\bY(t)=(\bM(t),\bQ(t))$
is a solution to (\ref{aveq}).
\smallskip\\
iii)  Let us  prove the  asymptotics
\be\la{intriateq}
 |\aM(t)-\bM^r|+|\aQ(t)-\bQ^r|\le C[p+d_0 e^{-p\mu t}],\qquad t\ge 0,
 \qquad{ p\le \ve},\quad{p}/\ga=r.
   \ee 
By the smoothness of the projection $\Pi$, $|\aY(0)-\bY^{r}|=|\Pi(\bX(0)-\bX^{r})|\le C d_0$ for small    $d_0= |X(0)-X^{r}|$.
 Moreover,  $\bY^r$ is a linearly stable  
 stationary  state of the averaged system (\ref{aveq}) by our assumptions. Therefore,
 extending 
Theorem 5.5.1 of \ci{SVM2007} onto 
the system (\ref{perN}) 
on the 
 phase space $\C\times\rS^2$ for   $\bZ(t)\in\rS^2$ and  $\aM(t)\in\C$,
 we obtain 
that for small $d_0$,
\be\la{intriateq2}
 \quad
 \sup_{t\in [0,\infty)}|\aY(t)-\bY(t)|=\cO(p),
 \qquad { p\to 0},\quad{p}/\ga=r.
   \ee 
  On the other hand, for small $d_0>0$ and  $\mu\in(0,\nu)$,   (\ref{lam12}) implies  that
    \be\la{Bt2}
 | \bY(t)- \bY^r|\le C d_0 e^{-p\mu  t},\qquad  t>0.
  \ee
 Hence,  (\ref{intriateq}) follows. Now (\ref{intriat}) holds by
(\ref{intr}).
 Finally, (\ref{intri}) is a particular case of (\ref{intriat}) when $d_0=0$.
    }

  \bc\la{cin}
  i)  Asymptotics (\ref{intriad}) hold for solutions $X(t)$ to the MBE with  $\Pi X(0)=(\bM_\pm,\bQ_\pm)$
 given by  (\ref{sol}) and  (\ref{nz22}).
   \smallskip\\
  ii)  Asymptotics (\ref{intriads}) hold for solutions with $\Pi X(0)=(\bM_\pm,\bQ_\pm)$, given by  (\ref{sol}), if the states are asymptotically stable for the  dynamics (\ref{aveq}).
   \smallskip\\
iii)  Asymptotics (\ref{intri}) and (\ref{intriat})  hold for solutions with  $\Pi X(0)=(\bM_+,\bQ_+)$ given by
 (\ref{nz22}).
  \ec
  
  \br\la{rape}\rm
  i) The asymptotic stability of  the harmonic states (\ref{sol}) 
  remains an open problem.
  \smallskip\\
  ii)
    For a suitable class of  {\it almost periodic pumpings} $A^e(t)$,
    asymptotics similar to (\ref{intriad})--(\ref{intriat})
  hold  if the corresponding order function
  $\de_\cD(p)=\cO(p^\nu)$ with $\nu\in(0,1)$, see  \ci[Section 4.6]{SVM2007}.
  \er

\section{Conflict of interest}
We have no conflict of interest.

 \section{Data availability statement}
The manuscript has no associated data.

  \appendix
  
  \protect\renewcommand{\thetheorem}{\Alph{section}.\arabic{theorem}}

\setcounter{equation}{0}
\section{The a priori bounds and well-posedness}\la{swp}
In this section, we 
prove the a priori bounds (\ref{apri})
 assuming that $A^e(t)\in C[0,\infty)$.
The bounds
imply
 the well-posedness of the MBE in
the phase space $\aX=\R^2\times \rS^3$. 
 The Schr\"odinger amplitudes $C_1(t), C_2(t)$
are bounded by the charge conservation (\ref{Blcc}).  Hence, it remains to prove 
the bounds for the Maxwell amplitudes $(A(t),B(t))$.   The following lemma 
and its proof refine the estimate (2.1) from \ci{K2024}.

\bl
There exists a Lyapunov function
$V(A,B)$ such that 
\be\la{aprL}
a_1[A^2+B^2]\le V(A,B)\le a_2[A^2+B^2]\quad {\rm where}\quad a_1,a_2>0,
\ee
 and for solutions to (\ref{HMB2}), the function $V(t)=V(A(t),B(t))$ satisfies
 the inequality 
\be\la{Vder3}
\dot V(t)\le -\ga b
V(t)
+d\fr{p^2}\ga,\qquad t>0;
\quad b,d>0.
\ee
 \el
\bpr
Denote by $E(A,B)=\fr1{2}(\Om^2A^2+B^2)$  ``the energy" of the Maxwell field.
The first two equations of (\ref{HMB2}) imply 
\be\la{denE}
\fr d{dt}E(A(t),B(t))=\Om^2A(t)B(t)
+B(t)[-\Om^2A(t)-\ga B(t)+cj(t)]
=
-\ga B^2(t)+ cB(t) {j(t)}.
\ee
We construct the Lyapunov function  following
the standard approach 
to dissipative perturbations of the Hamiltonian systems
 \ci{BV1992,H1981}:
$V(A,B)=E(A,B)+\ve AB$. Then (\ref{aprL}) holds for  small $\ve>0$. 

It remains to demonstrate (\ref{Vder3}) for sufficiently small $\ve>0$.
Differentiating, we obtain  
\beqn\nonumber
\dot V(t) &=&\dot E(t)+\ve\dot A(t)B(t)+\ve A(t)\dot B(t)\\
\nonumber
&=&-\ga B^2(t) +cB(t){j(t)} +\ve B^2(t)+\ve A(t)[-\Om^2 A(t)-\ga  B(t)+c j(t)]\\
\nonumber
&=&-(\ga-\ve) B^2(t)-\ve\Om^2 A^2(t)-\ve\ga A(t) B(t)+c{j(t)}(\ve A(t)+B(t)).
\eeqn
Note that $|\ve\ga AB|\le\fr\ga{2}B^2+\fr12\ga\ve^2 A^2$.
Hence,   for small $\ve>0$, 
\be\la{Vder2}
\dot V(t)\le-[\fr{\ga}{2}-\ve]B^2(t)-\fr\ve2 \Om^2A^2(t)+c{j(t)}(\ve A(t)+B(t)),\qquad t>0.
\ee
Moreover,  $|j(t)|\le \vka=p\om $ by (\ref{JOm}) and (\ref{Blcc}). 
Hence, for any $N>0$, 
$$c|j(t)(\ve A(t)+B(t))|\le \fr\ga N c^2(\ve A(t)+B(t))^2 + \om^2N\fr{ p^2}\ga.
$$
Therefore,  (\ref{Vder2}) with small $\ve<\ga/2$ and sufficiently large $N>0$ implies 
(\ref{Vder3}).
\epr

\bc
Solving the inequality (\ref{Vder3}),
 we obtain:
\be\la{Vdec}
V(t)\le 
V(0)+ \fr db r^2,\qquad t\ge 0,\quad r=p/\ga.
\ee
Hence,
for solutions to (\ref{HMB2}) with $p/\ga=r$,
 the following bounds hold:
 \be\la{apr2}
A^2(t)+B^2(t)\le D_r(A^2(0)+B^2(0)),\qquad t\ge 0.
\ee
Now (\ref{apri}) follows.
\ec

\setcounter{equation}{0}
  \section{On averaging theory on manifolds}\la{aC}
The averaging theory is usually considered for dynamical systems on linear spaces, 
 \ci{SVM2007}. On the other hand, definition  of the KBM vector field (\ref{KBM})
  and of the order function (\ref{KBM2}) 
 are invariant. Namely, let $\cM$ be a smooth finite-dimensional manifold,
 and consider the dynamical system
  \be\la{aMM}
 \dot Y(t)=\ve V(Y(t),t),\qquad t\ge 0,
 \ee
 where $V(Y)$ is a vector field on $\cM$. The corresponding
  averaged equation is defined by
 \be\la{aMMa}
  \dot \bY(t)=\ve\ov V(\bY(t)),\quad t\ge 0;\qquad \ov V(\bY)=\lim_{T\to\infty}\fr1T\int_0^TV(\bY,t)dt,\quad \bY\in\cM.
 \ee
 \bd i) (\ci[Definition 4.2.4]{SVM2007}) 
 $V(Y)$ is a KBM vector field if for any compact set $\cK\subset\cM$,
 \be\la{aKBM}
  \sup_{\bY\in \cK}\fr1{T}
 \Big| \int_0^T[V(\bY, t)-\ov V(\bY)]d t\Big|
  \to 0,\qquad T\to \infty.
  \ee  
  ii) The order function is defined as in \ci[Lemma 4.3.1]{SVM2007}:
     \be\la{aKBM2}
  \de_\cK(\ve):=\ve
  \sup_{\bY\in \cK}
  \,\,
  \sup_{T\in[0,\ve^{-1}]}
 \Big| \int_0^T[V(\bY, t)\!-\!\ov V(\bY)]d t\Big|.
  \ee

  \ed

  Accordingly, the 
  classical  Theorem 4.3.6 of \ci{SVM2007} admits the following extension.
Let the vector field $V(Y,t)$ be Lipschitz continuous uniformly in time: 
in each local chart on every subset $K\subset \cM$,
\be\la{LIP}
|V(Y_1,t)-V(Y_2,t)|\le L(K)|Y_1-Y_2|, \qquad Y_1,Y_2\in K,\quad t\ge 0,
\ee
where $L(K)<\infty $ can depend also on the choice of the local coordinates.
Note that (\ref{LIP}) implies (\ref{aKBM}).

\bt\la{tC} 
Let $\cM$ be a compact manifold,  (\ref{LIP}) hold, $Y(t)=Y_\ve(t)$ be solutions to (\ref{aMM})
with $\ve\in(0,1]$ and the same initial state $\bY_\ve(0)=Y^0$, 
and $\bY(t)=\bY_\ve(t)$ be solutions to (\ref{aMMa}) with $\bY_\ve(0)=Y^0$.
Then
\be\la{adiab}
\max_{t\in[0,\ve^{-1}]}\rho(Y_\ve(t),\bY(t))=\cO(\de^{1/2}),\qquad\ve\to 0,
\ee
where $\rho$ is an arbitrary distance defined on $\cM$.
\et
\bpr The compact manifold $\cM$ is 
 diffeomorphic to a smooth submanifold of $\R^N$ with a suitable $N\ge 1$.
Now the bounds (\ref{adiab}) follow by the same arguments as Theorem 4.3.6 of \ci{SVM2007} 
with the interpretation of points $Y\in\cM$ and tangent vectors $V(Y,t)\in T_Y\cM$
as the corresponding  vectors from $\R^N$. 
With this interpretation,
the vector field $V(Y,t)$ on the compact submanifold $\cM$
is Lipschitz continuous by   (\ref{LIP}).
\epr

\setcounter{equation}{0}
\section{The Maxwell--Bloch equations for pure states}\la{sMB}
 In \ci{K2024}, the 
 MBE  is obtained as the Galerkin approximation of the coupled
Maxwell--Schr\"odinger system. 
The approximation
consists of a single-mode Maxwell field coupled to
  two-level molecule in a bounded cavity $V\subset\R^3$:
 \be\la{solMB2}
\bA(x,t)=A(t)\bX(x),\quad
\psi(x,t)=C_{1}(t)\psi_{1}(x)+C_{2}(t)\psi_{2}(x),\qquad x\in V.
\ee
 Here 
 $\bA(x,t)$ denotes the vector potential of the Maxwell field, and
 $\bX(x)$ is a normalised eigenfunction of the Laplace
 operator in $V$ under suitable boundary value conditions
 with an eigenvalue  $-\Om^2/c^2$. By
$\psi_{l}$ we denote  some 
normalised eigenfunctions of the 
Schr\"odinger operator $\bH:=-\fr{\hbar^2}{2\cm} \Delta  + e\Phi (x)$ with the corresponding 
eigenvalues $\hbar\om_1<\hbar\om_2$, where $\Phi (x)$
is the molecular (ion's) potential.
The MBE  read as the Hamiltonian system with a dissipation and an external source:
\be\la{HMB} 
\fr1{c^2}\dot A(t)=\pa_B H,\quad \fr1{c^2}\dot B(t) =-\pa_A H-\fr\ga{c^2} B;
\qquad i\hbar\dot C_{l}(t)=\pa_{\ov C_{l}}H,\quad l=1,2.
\ee
The Hamiltonian $H$ is defined as
$$
H(A,B,C_1,C_2,t)=\cH(A\bX,
B\bX, C_{1}\vp_{1}
+C_{2}\vp_{2},t),
$$
where
$\cH$ is the Hamiltonian of the coupled
Maxwell--Schr\"odinger equations 
with pumping \ci{BT2009,GNS1995,K2019phys,K2013,K2022,KK2020,NW2007,S2006}.
Neglecting the spin and scalar potential (which can be easily added),
the Hamiltonian $H$,  in the traditional {\it dipole approximation},
reads as \ci[(A.5)]{K2024}:
\be\la{Hc33}
H(A,B,C_1,C_2,t)=\fr1{2c^2}[B^2
+ \Om^2A^2]+\hbar\om_1 |C_{1}|^2+
\hbar\om_2|C_{2}|^2
  -\fr {2\vka}c
[A+A^e(t)]\,
 \rIm[\ov C_{1}C_{2}].
\ee
 Now the Hamilton equations (\ref{HMB}) become (\ref{HMB2}).
 The charge conservation 
 $|C_{1}(t)|^2+|C_{2}(t)|^2=\const$
 follows by differentiation from the last two equations of  (\ref{HMB2})
since the function $a(t)$ is real-valued. 
We  consider solutions with $\const=1$, see (\ref{Blcc}).

\setcounter{equation}{0}
  \section{On possible treatment of the  laser action}\la{aD}
  Here we discuss possible treatment of the laser action 
  relying on the obtained results.
  \smallskip\\
  {\bf On the smallness of the parameters.}
Note that  the dipole moment $p$ and dissipation coefficient $\ga$ are very small 
for many types of lasers. In particular, 
the dissipation coefficient  for the Ruby laser is the electrical conduction of corundum
which is 
 $\ga\sim10^{-14}$ in the Heaviside--Lorentz 
 units
\ci{DH1992,S1986,S2012}.
For the dipole moment typically 
$|p|\sim 10^{-18}$ according to
\ci{MW}
that agrees with the classical dipole moment $ed/2$, where 
$d\sim 10^{-8}$cm is the molecular diameter, and
$e\sim 10^{-10}$ is the elementary 
charge in the same units. 
  \smallskip\\
{\bf On the laser threshold.}
The asymptotics (\ref{intriad}) and (\ref{intri}) hold for solutions with the harmonic initial states $(\bM^r_\pm,\bS^r_\pm)$. On the other hand,
(\ref{intriads}) and (\ref{intriat}) hold for solutions with  initial states 
from an open domain of attraction in the phase space. 
Hence, 
the asymptotics (\ref{intriads}) and (\ref{intriat}) appear with a ``nonzero probability" in contrast to 
 (\ref{intriad}) and (\ref{intri}). This fact, provisionally, 
 shed new light on the ``laser threshold" which is necessary
 to ignite the laser action: the  intensity of random pumping must be sufficiently large to bring the solution to
  the domain
 of attraction, and then the solution is captured in the domain with the single-frequency asymptotics.
\smallskip\\
{\bf On the laser amplification.}
The equations (\ref{HMB2}) describe  one molecule  coupled to
 the Maxwell field. 
The limiting amplitudes of the Maxwell field
 in all asymptotics 
(\ref{intriad})--(\ref{intri})
   do not depend on non-resonance harmonics in the pumping (\ref{qpp})
  with the frequencies $\Om_k\ne\om$.
 This means that the dynamics  (\ref{HMB2})  acts as a filter, selecting only the resonant harmonics, that itself cannot explain the
  amplification of the Maxwell field in laser devices.
 The amplification could be explained by a large number of active  molecules,
 typically  $N\sim 10^{20}$, under the traditional assumption that the molecules
 interact with the Maxwell field but do not interact with each other \ci{N1973}.
 In this case, the amplitude of the total Maxwell field is multiplied by $\sqrt{N}\sim 10^{10}$ by the Law of Large Numbers. 
The amplification is possible only due to the single-frequency asymptotics of individual contributions.

The amplification  also  depends significantly on the fact that  
the phases of the amplitudes  $\bM^r$
 in  asymptotics  (\ref{intriad})--(\ref{intriat})
for all molecules are  uniformly distributed. 
 The uniformity   takes place
 if it holds for phases of  $\bA^e$.
This follows from the asymptotics and 
Remark \ref{runif} since $e^{i\phi}e^{-i\Om t}\bM$ differs from $e^{-i\Om t}\bM$
by a shift of time.
\smallskip\\
{\bf Self-induced transparency.} The harmonic states (\ref{S1})  means that the incident wave 
$\bA^e e^{-i\Om t}$ results in the outgoing wave $M(t)\approx -\bA^e e^{-i\Om t}$ with the same 
amplitude and frequency, but with the phase jump $\pi$. This phenomenon resembles 
 the self-induced transparency \ci{S1986}.

 \end{document}